\documentclass{article}
\usepackage{graphicx}
\usepackage{amsmath,amssymb,amsthm,mathtools}
\usepackage[colorlinks=true]{hyperref}
\usepackage{makeidx}
\usepackage[symbol]{footmisc}
\usepackage[utf8]{inputenc}

\usepackage{subfig}
\usepackage{float}
\usepackage{multirow}
\usepackage{longtable,lscape}
\usepackage{cases}
\theoremstyle{plain}
\newtheorem{theorem}{Theorem}
\newtheorem{remark}[]{Remark}
\usepackage{authblk}

\usepackage[nice]{nicefrac}

\title{Multi-fidelity and multi-level Monte Carlo methods for kinetic models of traffic flow}
\author{Elisa Iacomini\footnote{Department of Environmental and Prevention Sciences, University of
		Ferrara, Italy, \textit{elisa.iacomini@unife.it}.}  and Lorenzo Pareschi\footnote{Maxwell Institute for Mathematical Sciences and Department of Mathematics, School of Mathematical
		and Computer Sciences (MACS), HWU Edinburgh, UK.\\
		Department of Mathematics and Computer Science, University of Ferrara, Italy, \textit{lorenzo.pareschi@unife.it}.}}%{E.~Iacomini, L.~Pareschi }
\begin{document}
\maketitle

%------
% Insert the title of your paper and (if necessary)
% a short title for the running head.
%------

%------
% Insert full names of the authors.
% Add further authors as follows:
%  \emsauthor{2}{}{}
%  \emsauthor{3}{}{}
% etc.
% Abbreviate first names for the running head.
%------

%------
% Use \authormark if the list of authors is too
% long for the running head: \authormark{A.~Doe et al.}
%------

%------
% Add one \emsaffil and one \email for each author.
% NOTE: The address does NOT appear in the paper.
% It will probably be printed in an appendix.
%------

%------
% Add MSC 2020 codes according to www.ams.org/msc/msc2020.html.
% Secondary codes (in square brackets) are optional.
%------

%------
% Add a list of keywords.
%------

%------
% Optional: dedication
%------
%\chapterdedication{Dedicated to ...}

%------
% Insert your abstract.
%------
\begin{abstract}
In traffic flow modeling, incorporating uncertainty is crucial for accurately capturing the complexities of real-world scenarios. 
In this work we focus on kinetic models of traffic flow, where a key step is to design effective numerical tools for analyzing uncertainties in vehicles interactions. To this end we discuss space-homogeneous Boltzmann-type equations, employing a non intrusive Monte Carlo approach both on the physical space, to solve the kinetic equation, and on the stochastic space, to investigate the uncertainty.
To address the high dimensional challenges posed by this coupling, 
control variate approaches such as multi-fidelity and multi-level Monte Carlo methods are particularly effective. While both methods leverage models of varying accuracy to reduce computational demands, multi-fidelity methods exploit differences in model fidelity, while multi-level methods utilize a hierarchy of discretizations. Numerical simulations indicate that these approaches provide substantial accuracy improvements over standard Monte Carlo methods. Moreover, by using appropriate low-fidelity surrogates based on approximated steady state solutions or simplified BGK interactions, multi-fidelity methods can outperform multilevel Monte Carlo methods. 
\end{abstract}


%\makecontribtitle

%------
% INSERT THE BODY OF THE PAPER HERE (except
% acknowledgments, funding info and bibliography)
%------
\tableofcontents

\section{Introduction}

Traffic flow modeling is a critical component in the design, management, and optimization of transportation systems. One of the most powerful mathematical tools used in this field is partial differential equations (PDEs). PDEs offer a robust framework for describing the dynamics of traffic flow, capturing essential features such as density, speed, and flux over time and space. These equations enable researchers to predict and analyze complex traffic patterns, facilitating the development of effective traffic management strategies and infrastructure improvements. In this setting, kinetic traffic flow models provide a statistical description of traffic \cite{klar2000kinetic,herty2020bgk, puppo2017kinetic,albisurv19}, taking into account both car-to-car interactions and mass distribution of traffic. 

However, real-world traffic is inherently uncertain due to various factors such as fluctuating demand, unpredictable incidents, and diverse driver behaviors. This uncertainty can significantly impact the accuracy and reliability of traffic flow models \cite{herty2021reconstruction,herty2018hybrid,puppo2016fundamental,tosin2021uncertainty}. Therefore, incorporating uncertainty into PDE-based traffic models is crucial for creating more realistic and trustful solutions. 

Several approaches to quantify uncertainty are presented in the literature and can be classified in non-intrusive and intrusive methods. The main idea underlying the former approach is to solve the model for fixed number of samples using deterministic numerical algorithms. Typical examples are Monte Carlo (MC) and stochastic collocation methods \cite{bertaglia2022bi,tosin2021uncertainty,pareschi2021,dimasurv24}.
On the other side, intrusive approaches are based on the fact that the governing equations have to be modified to incorporate the probabilistic character of the model parameters \cite{S21,Gottlieb2001}. The stochastic Galerkin method is one of the most famous in this framework. Here, stochastic processes are represented as piecewise orthogonal functions, also known as generalized polynomial chaos expansion (gPC) \cite{S2,S16, S1,S3}, and then substituted into the governing equations. A Galerkin projection is then used to obtain deterministic evolution equations for the coefficients of the series expansions \cite{pettersson2014stochastic}. 

This approach has recently been applied to traffic flow models, with studies investigating its {performance} across various scales of observation \cite{iacomini2023overview}. In particular the uncertainty is introduced at microscopic, mesoscopic and macroscopic scales, and the resulting stochastic models are analyzed. Many challenges arise here, since some desired properties of the original system are not necessarily transferred to the intrusive formulation, in particular at the macroscopic level we face the loss of hyperbolicity of the system \cite{jin2019study,poette2009uncertainty}.  To preserve hyperbolicity, the basis functions must meet additional assumptions, and a consistent gPC expansion is required \cite{gerster2021stability}. Additionally, by linking with the kinetic model, the probability of high-risk traffic zones where instabilities may occur can be studied \cite{herty2022uncertainty}.

However, in many practical scenarios the uncertainty distribution is either unknown or irregular, posing challenges for the stochastic Galerkin method, which relies on regularity.  As a result, non-intrusive methods are often better suited for these applications. These methods not only preserve the structure of the underlying numerical solver but are also easier to parallelize \cite{xiu2010}, while mitigating the curse of dimensionality typically encountered in uncertainty analysis. This is particularly relevant when one deals with kinetic equations of traffic flow with the general form
\begin{equation}\label{eq:kin}
	\partial_t f + v \cdot \partial_x f =  \frac1{\tau} Q(f,f)
\end{equation}
where $f=f(t,x,v;z)$ is the distribution function, at time $t\ge0$, space $x \in \mathcal{D}\subseteq \mathbb{R}$, velocity $v \in \mathbb{R}$, and $z\in \Omega \subseteq \mathbb{R}^{d_z}$, $d_z \ge 1$ is a random variable. The parameter $\tau \geq 0$ characterizes the traffic flow regime and plays the role of the Knudsen number. The particular structure of the interaction term $Q(f,f)$ depends on the kinetic model considered and will be discussed in the Section 2. 

Despite the advantages, the common prototype for non-intrusive approach, i.e. the Monte Carlo method, suffers from slow convergence. To overcome this limitation, several strategies have been developed \cite{niederreiter1992random,caflisch1998monte}. Among those, multi-fidelity methods based on multi-scale models \cite{dimarco2019multi} and the multi-level Monte Carlo approach \cite{giles2015multilevel} have shown to be the most effective in terms of flexibility and generality.
The former approach exploits the multi-scale nature of the problem to effectively reduce the  variance in Monte Carlo simulations along the different scales. Specifically, we exploit the hierarchical relationship between different model: high-fidelity models, such as kinetic models,  offer high accuracy but are computationally expensive, while low-fidelity models, such as macroscopic models, are less precise but computationally less demanding. The crucial aspect is that the low fidelity model must maintain correlation with the high fidelity model in the space of uncertainties. By performing a limited number of high-fidelity evaluations and numerous low-fidelity ones, we can improve accuracy without significantly increasing computational costs.
Instead of performing all simulations at a single resolution, multi-level Monte Carlo combines simulations at multiple levels of accuracy, where coarse simulations are inexpensive and provide broad trends, while fine simulations are used to correct errors. The key insight is that by carefully balancing the number of simulations across different levels, the computational cost of achieving a given accuracy can be significantly reduced compared to standard Monte Carlo methods. This makes this approach especially powerful in high-dimensional problems, as in the uncertainty quantification framework.

An important aspect of this work is that our methodology is based on a Direct Simulation Monte Carlo (DSMC) solver \cite{pareschi2022mean}, rather than a deterministic solver for the kinetic model \cite{dimarco2019multi}. This choice is motivated by the fact that Monte Carlo solvers are widely used in order to solve Boltzmann type equations and can be easily applied to different models by exploiting the rules that define the microscopic dynamics. However, this approach introduces the additional challenge of coupling the DSMC solution of the traffic model with the Monte Carlo solver to address uncertainty. Additionally, given the popularity of splitting algorithms in DSMC methods that separately address the free motion of vehicles and their interactions, we focus on space homogeneous models. In these models, we disregard spatial dependency, treating the system as uniform across space. This simplification enables us to isolate and analyze the core dynamics of vehicle interactions, which are governed by the vehicle interaction operator—the most computationally expensive part of the problem, without the added complexity introduced by spatial variation.

The rest of the paper is organized as follows: Section\ref{sec:model}  introduces the kinetic traffic model. The variance reduction approaches are presented in Section \ref{sec:CV}. Specifically, in Subsection \ref{sec:multifidelity}, we outline the main principles of the multi-fidelity approach, while Subsection \ref{sec:multilevel} covers the multi-level Monte Carlo method. Finally, Section \ref{sec:appl} presents the application of these approaches to the traffic model and discussed the results of the numerical experiments.

\section{Kinetic models of traffic with uncertainties}\label{sec:model}
%\subsection{Homogeneous case}

A kinetic description of traffic flow is based on the identification of the interaction rules. In this work we assume the interaction rules as in \cite{tosin2021uncertainty}, in particular, we consider pairwise interactions which modify the speed of the vehicles involved. 

Let us denote by $v\in [0,1]$ the (renormalized) speed of a vehicle and by $v_*\in [0,1]$ the (renormalized) speed of the leading vehicle, where for leading vehicle we mean the vehicle in front. The update of the velocities due to the interaction is described by:
\begin{align}
	\label{eq:update1}
	&v'=v+\gamma I(v,v^*;z) + \sigma D(v)\eta \\
	\label{eq:update2}
	&v'_*=v_*  
\end{align}
where $\gamma>0$ is a proportionality parameter and $I$ is the interaction function which depends on the speeds before the interactions and on an uncertain parameter, a random variable $z\in \mathbb{R}$ with known probability distribution. 
In particular $\eta$ is a random variable with $0$ mean and unitary variance, $\sigma>0$ characterizes the intensity of the noise, and $D(\rho,v)=a(\rho)\sqrt{v(1-v)}$, $a(\rho)=\rho(1-\rho)$ models the diffusion process caused by the intrinsic randomness of the driver behaviour.  %Please be aware that while $z$ is the stochastic variable of the model, $\eta$ describes the noise of the interactions.
Accordingly to \cite{tosin2019kinetic}, we consider the following interaction rule:
\begin{equation} \label{eq:interaction}
	I(v,v_*;z)= P(\rho;z)(1-v) + (1-P(\rho;z))(P(\rho;z)v_*-v)
\end{equation}
where $\rho \in [0,1]$ is the density of the vehicles, $P(\rho;z)\in [0,1]$ is the probability of acceleration, as in \cite{tosin2021uncertainty}, which depends on the traffic density, the higher $\rho$ is, the lower $P$, since in heavy traffic conditions preclude accelerations, and vice versa.
 $P$ is also uncertain, since it depends on $z$:
\begin{equation}
	P(\rho;z)=(1-\rho)^z
\end{equation}
where we assume $z>0$.

We refer to \cite{tosin2021uncertainty} for details on the admissibility of the interaction rules \eqref{eq:update1}-\eqref{eq:update2}. Here we just recall the sufficient conditions:
\begin{equation} \label{eq:suff}
|\eta|\le c(1-\gamma), \qquad cD(v)\le \text{min}\{v,1-v\}
\end{equation}
where $c>0$ is an arbitrary constant.

In order to describe the aggregate dynamics given by the superposition of the many binary interactions, we introduce the distribution function $f(t,v;z)$, such that $f(t,v;z)dv$ describes the probability of a vehicle to travel with a speed spanning between $v$ and $v+dv$ at time $t$ given the uncertain parameter $z$. 

As stated in \cite{tosin2021uncertainty}, the distribution function evolves as a space homogeneous Boltzmann-type equation for Maxwellian-like particles where the interaction probability is uniform. In its weak form, this equation is expressed as
\begin{equation}\label{eq:Boltzmann}
\begin{split}
	\frac{d}{dt}\int_0^1 \varphi(v)f(t,v;z)dv&=\frac{1}{\tau}\int_0^1 \varphi(v)Q(f,f)(t,v;z)dv\\
	&=\frac{1}{2\tau}\int_0^1\int_0^1
\langle\varphi(v')-\varphi(v)\rangle f(t,v;z)f(t,v_*;z) dv_*dv
\end{split}
\end{equation}
for every observable quantity $\varphi:[0,1]\to\mathbb{R}$ which refers to any quantity that can be represented as a function of the microscopic state $v$ of the vehicle. Note that \eqref{eq:Boltzmann} is a stochastic kinetic equation, since $f$ depends on the uncertain parameter $z$. For more details on \eqref{eq:Boltzmann} we refer to \cite{tosin2021uncertainty}.

\subsection{Approximated steady states}

We are particularly interested to the asymptotic behavior and an explicit computation of the steady states. This is challenging for the Boltzmann type equation just described. For this reason we approximate the Boltzmann-type equation with a Fokker Planck partial differential equation through an analogous of the classical grazing collisions limit \cite{pareschitoscani2013,hertypar10}. Let us assume $\gamma,\sigma^2<<1$, namely the deterministic part of the interactions and the stochastic fluctuations are small. On the other hand, the frequency of the interactions has to increase accordingly, so that we obtain the following scaling $\gamma=\varepsilon$, $\sigma^2=\lambda\varepsilon$ and $\tau={\varepsilon}/{2}$, where $\varepsilon$ is the scaling parameter and $\lambda>0$ a constant of proportionality.
 
In this way, we have the scaled interaction rules 
\begin{align}
	\label{eq:update1_scaled}
	&v'=v+\varepsilon I(v,v^*;z) + \sqrt{\varepsilon  \lambda} D(v)\eta \\
	\label{eq:update2_scaled}
	&v'_*=v_* .
\end{align} 

Following the computation in \cite{tosin2021uncertainty}, in the limit $\varepsilon\to 0$ the traffic flow dynamic is well-described by the Fokker-Planck equation  
\begin{equation}
	\partial_t f=\frac{\lambda}{2}\partial^2_v(D^2(v)f)-\partial_v\left[(P(1+(1-P)U(t;z))-v)f\right]
\end{equation}
with
\[
U(t;z) = \int_{0}^{1} v f(t,v;z)\,dv.
\]
 Direct computations show that the corresponding steady state is given by
 \begin{equation}\label{eq:steady}
 	f^\infty(v)= \frac{v^{\frac{2U^\infty (\rho;\xi)}{\alpha}-1}(1-v)^{\frac{2(1-U^\infty)}{\alpha}-1}}{\textup{Beta}\left( \frac{2 U^\infty}{\alpha}, \frac{2(1-U^\infty)}{\alpha}\right)},
 \end{equation}
where $\alpha=\lambda a^2(\rho)$, $\textup{Beta}$ stands for the Beta distribution, and
\begin{equation}\label{eq:eq_vel}
	U^\infty(\rho)=\frac{P(\rho;z)}{P(\rho;z)+(1-P(\rho;z))^2}.
\end{equation}

The analytical knowledge of the steady state will be crucial in what follows, since it will represent the key element in the design of low-fidelity surrogates.

\section{Multi-fildelity and multi-level methods}\label{sec:CV}
Before going into the details, we introduce some notations that will be used throughout the paper.

If $z\in \Omega$ is distributed as $p(z)$, the expected value, or mean, of $f(t,v;z)$ will be denoted by
\begin{equation} 
	\mathbb{E}[f](t,v)=\int_\Omega f(t,v;z)p(z)dz, 
\end{equation}
and the variance of $f(t,v;z)$ is defined as 
\begin{equation}
	\textup{Var}(f)(t,v) = \int_\Omega \left(f(t,v;z)-\mathbb{E}[f](t,v)\right)^2 p(z) dz.
\end{equation}

In the physical space we adopt a Direct Simulation Monte Carlo (DSMC) approach to solve the kinetic equation for traffic flow. Here we do not describe the details of the method and we refer to \cite{KHP2005,pareschitoscani2013,pareschi2022mean} for a more in depth presentation. 

We briefly recall that for a DSMC simulation using $N$ samples $\{v_1,\ldots,v_N\}\in [0,1]$ at time $t>0$ we denote with $f_{N}$ the histogram reconstruction with uncertainty as
\begin{equation}\label{eq:reconstr}
	f_{N,\Delta v}(t,v;z)=\frac{1}{N}\sum_{i=1}^N S_{\Delta v} (v -v_i(t;z)),
\end{equation} 
where $\Delta v$ is the mesh size and $S_{\Delta v} (\cdot)$ a suitable approximation of the Dirac delta function $\delta(\cdot)$ characterizing the empirical measure.
In the simplest case, we have $S_{\Delta v}(v)=\chi(|v|\leq \Delta v/2)/\Delta v$, where $\chi(\cdot)$ is the indicator function, that corresponds to the standard histogram reconstruction.

Note that, uncertain observables of the kinetic distribution 
\[
(\varphi, f)(t;z) = \int_0^1 \varphi(v)f(t,v;z)dv, 
\]
do not need any reconstruction and can be evaluated directly from the empirical measure $f_N(t,v;z)$ as
\[
(\varphi,f_N) (t;z) = \dfrac{1}{N} \sum_{i = 1}^N \varphi(v_i(t;z)). 
\]
%Hence, by assuming that $\int_V f(t,w)dw= 1$ we have that $(\varphi,f)=\mathbb E_V[\phi]$, where $\mathbb E_V[\cdot]$ is the expectation of the observable quantity $\varphi$ with respect to the density $f$. In the sequel, we will also make use of the notation $\mathbb E[\cdot]$ to denote the expectation in the random space of uncertainties. We will implicitly assume that for multidimensional variables, as in the case of statistical samples, the expected values are done with respect to each variable.

%Thanks to the central limit theorem the following result holds \cite{caflisch1998monte}:
%\begin{lemma}\label{lem:1}
%The root mean square error is such that for each $t\ge 0$
%\begin{equation*}
%\label{eq:th_e1}
%\mathbb E_{V}\left[\left( (\varphi,f) - (\varphi,f_N) \right)^2\right]^{1/2}= \dfrac{\sigma_\varphi}{N^{1/2}}, 
%\end{equation*}
%where $\sigma^2_\varphi = \VV_V[\varphi]$ with
%\begin{equation*}
%\label{eq:sigma_varphi}
%\VV_V[\varphi](t) = \int_V( \varphi(w)- (\varphi,f)(t))^2 f(t,w)dw.
%\end{equation*}
%\end{lemma}

The numerical error of the reconstructed DSMC solution \eqref{eq:reconstr}, by ignoring the time discretization error and in the in the case where the kinetic density does not depend on the uncertainties $f=f(t,v)$, can be estimated from \cite{pareschi2022mean}
\begin{theorem} The error introduced by the reconstruction function \eqref{eq:reconstr} satisfies
\begin{equation}
\left\| f(t,\cdot)-f_{N,\Delta v}(t,\cdot)\right\|_{L^p([0,1],L^2([0,1]))} \leq \frac{C_{S,f}}{N^{1/2}} + C_f (\Delta v)^q,
\label{eq:est1}
\end{equation}
accordingly to the order of accuracy $q$ used in the histogram reconstruction. In the above estimate, $C_f$ depends on the $q$ derivative in $v$ of $f(t,v)$ and $C_{S,f}$ depends on $S_{\Delta v} (\cdot)$ and $f$.
\label{th:1}
\end{theorem}
In \eqref{eq:est1} for a function $g(v;v_1,\ldots,v_N)$ we used notation
\begin{equation}
\|g\|_{L^p([0,1],L^2([0,1]))}=\|\mathbb E_{v}\left[g^2\right]^{1/2}\|_{L^p([0,1])},
\label{eq:norm1}
\end{equation}
where $\mathbb E_{v}$ is the $N$-dimensional expectation of $g$ with respect to $v_i$, $i=1,\ldots,N$ identically distributed as $f(t,v)$.

\subsection{Standard MC sampling}
With the above notation the standard MC sampling steps for uncertainty quantification consists in the following three steps.

\begin{itemize}
	
	\item[] 1 - {\it Sampling}: Sample $M$ independent identically distributed (i.i.d.) values of $z$ from the distribution $p(z)$
	\item[] 2 - {\it Solving}: For each realization of $z_k$, $k=1,\dots,M$ the underlying kinetic model is solved numerically with the DSMC solver. We denote the histogram reconstruction on the mesh $\Delta v$ at time $t^n$ by $f^{n}_{k,N}(v)=f_{N,\Delta v}(t^n,v;z_k)$, $k=1,\dots,M$ where $N$ is the sample size of the DSMC solver. 
	\item[] 3 - {\it Reconstruction}: Estimate the expectation of the quantity of interest, in our case the kinetic density $f^n$, of the random solution:
	\begin{equation}
		\mathbb{E}[f^n] \approx E_M[f^n_{N,\Delta v}] \coloneq\frac{1}{M}\sum_{k=1}^{M} f^n_{k,N}.
		\label{eq:estimate}
	\end{equation}
\end{itemize}

%and the empirical mean, or Monte Carlo estimator, will be represented as
%\begin{equation} 
%	E_M[f]=\frac{1}{M}\sum_{k=1}^M f^k =\frac{1}{M}\sum_{k=1}^M f(t,v,z^k). 
%\end{equation}
%where $M$ is the number of samples.

To analyze uncertainty, the standard non-intrusive method is the Monte Carlo technique, which is known to suffer from slow convergence \cite{loeve1977probability}. 
In fact, we recall the following classical result \cite{caflisch1998monte}.

\begin{theorem}
	The root mean square error is such that for each $t\ge0$
	\begin{equation} \label{eq:MC_rate}
		\mathbb{E}\left[(\mathbb{E}[f]-E_M[f])^2\right]^\frac{1}{2}=  \frac{\sigma_f}{M^{\frac{1}{2}}}
	\end{equation}
	where $\sigma_f^2=\textup{Var}(f)$
\end{theorem}

The following theorem is a direct consequence of Theorem \ref{th:1} and the above result \cite{pareschi2022mean}.

%\begin{theorem} 
%	The root mean square error satisfies 
%	\begin{equation}
%		\mathbb{E}\left[ \mathbb{E}_V \left[ (\mathbb{E}[f] - E_M[f_N])^2\right] \right]^\frac{1}{2} \le \frac{\sigma_f}{M^\frac{1}{2}} + \frac{\nu_M}{N^\frac{1}{2}}
%	\end{equation}
%	where $\sigma^2 = \textup{Var}(f)$ and $\nu^2_M=E_M[\sigma^2_v]$ with $\sigma^2_v$ the variance with respect to the kinetic equation.
%\end{theorem}
%
%The proof can be found in \cite{pareschi2022mean}.

%Let us focus on the accuracy of the expectation of the solution $\mathbb{E}[f]$. 

\begin{theorem}\label{thm:rec}
	The error introduced by the MC estimate \eqref{eq:estimate} satisfies
	\begin{equation}
		|| \mathbb{E}[f](t, \cdot) - E_M[f_N](t, \cdot)||_{L^p([0,1],L^2(\Omega))} \le \frac{\tilde C_{S,f}}{M^{\nicefrac{1}{2}}} + \frac{\tilde C_{S,f,M}}{N^{\nicefrac{1}{2}}} + C_{{\mathbb E}[f]} (\Delta v)^q
	\end{equation}
	where $C_{{\mathbb E}[f]}$ depends on the $q$ derivative in $v$ of ${\mathbb E}[f]$, $\tilde C_{S,f}$ depends on $S_{\Delta v} (\cdot)$ and $f$, and $\tilde C_{S,f,M}$ depends on $S_{\Delta v} (\cdot)$, $f$ and $M$.
\end{theorem}
\noindent
%In case of histogram reconstruction the last term can be replaced by $C \Delta v^q$.

It is evident that the variance plays a crucial role in the error estimate, meaning that applying techniques to reduce the variance will improve convergence.

\subsection{Multi-fidelity approach} \label{sec:multifidelity}

To explore the uncertainty $z$ present in the traffic flow model, we first review the key concepts of the control variate multi-fidelity approach as introduced in \cite{dimarco2019multi,dimarco2020multiscale}.

%% THIS HOLDS IF THE Q(f,f) IS DETERMINISTIC AND THE INITIAL DATA IS RANDOM

%Moreover, it has been proved in \cite{dimarco2019multi} that the Monte Carlo estimator satisfies the following error bound in the $L^1$ norm with polynomial weight in the physical space:
%\begin{equation}
%	|| \mathbb{E}[f](\cdot)
%\end{equation}

%To improve the performance of the standard Monte Carlo method, one approach is to reduce the variance in \eqref{eq:MC_rate}, employing variance reduction techniques. 
The main idea is to exploit an approximate solution that shares a similar asymptotic behavior as the original, such as an approximated steady state of the kinetic equation, to accelerate convergence. Observe that the steady state \eqref{eq:steady} corresponds to the steady state of the Fokker Planck equation, which converges to the steady state of the Bolzmann equation only in the limit for $\varepsilon \to 0$. For this reason, we denote the solution of the time scaled Boltzmann equation by $f_{\varepsilon}(t,v;z)$ and the solution of the corresponding Fokker Planck equation as $\tilde{f}=\tilde{f}(t,v;z)$, as in \cite{pareschi2022mean}. The time scaled kinetic equation for small values of the scaling parameter $\varepsilon$ satisfies:
\begin{equation}
\lim_{\varepsilon \to 0} f_{\varepsilon}(t,v;z) =\tilde{f}(t,v;z).
\end{equation}
Therefore
\begin{equation}
\lim_{t \to \infty}	\lim_{\varepsilon \to 0} f_{\varepsilon}(t,v;z) =\tilde{f}^{\infty}(v;z).
\end{equation}

This step is crucial because the steady state $\tilde{f}^{\infty}(v;z)$ can often be computed analytically, as in this case \eqref{eq:steady}, while the steady state of the Boltzmann equation \eqref{eq:Boltzmann} remains unknown.

If we denote by $q[\cdot]$ any  \textit{quantity of interest}, the parameter dependent control variate approach with $\lambda \in \mathbb{R}$ can be formulated as 
\begin{equation}\label{eq:control_variate}
 q^\lambda[f_\varepsilon] =  q [f_\varepsilon] - \lambda \left( q[\tilde{f}] - \mathbb{E}[q[\tilde{f}]]\right).
\end{equation}
%The notation $q$ is used to emphasize the general applicability of the method to any  \textit{quantity of interest}. 
It is clear that the expected value fulfills
\begin{equation}
	\mathbb{E}[q^\lambda[f_\varepsilon]] = \mathbb{E}[q[\tilde{f}]].
\end{equation}

Then, we recall the following Theorem from \cite{pareschi2022mean}.

\begin{theorem}\label{thm:multi}
	If $\textup{Var}(q[\tilde{f}])\neq 0$, then the quantity 
	\begin{equation}\label{eq:lambda_star}
		\lambda^*=\frac{\textup{Cov}(q[f_\varepsilon],q[\tilde{f}])}{\textup{Var}(q[\tilde{f}])}
	\end{equation}
	minimizes the variance of $q^\lambda[f_\varepsilon]$ at the point $(t,v)$ and gives
	\begin{equation}\label{eq:var_lambdastar}
		\textup{Var}(q^{\lambda^*}[f_\varepsilon])=(1-\rho^2_{q[f_\varepsilon],q[\tilde{f}]})\textup{Var}(q[f_\varepsilon])
	\end{equation}
	where $\rho_{q[f_\varepsilon],q[\tilde{f}]}\in [-1,1] $ is the correlation coefficient of $q[f_\varepsilon]$ and $q[\tilde{f}]$. In addition we have 
	$$
	\lim_{\varepsilon\to 0} \lambda^*(t,v)=1, \qquad  \lim_{\varepsilon\to 0} \textup{Var}(q^{\lambda^*}[f_\varepsilon])=0 \quad \forall v \in \mathbb{R}^{d_v}.
	$$
	
\end{theorem}

Note that \eqref{eq:var_lambdastar} may lead to a strong variance reduction  when the correlation coefficient approaches one.

In the easiest case, the control variate formulation \eqref{eq:control_variate} can be modified using the steady state $\tilde{f}^\infty(v,z)$ of the Fokker-Planck equation:
\begin{equation}\label{eq:cv_steady}
	q^\lambda[f_\varepsilon]=q[{f}_\varepsilon] - \lambda \left( q[\tilde{f}^\infty] - \mathbb{E}[q[\tilde{f}^\infty]]\right).
\end{equation}

Here,
 $$
 \lambda^* = \frac{\textup{Cov}(q[f],q[\tilde{f}^\infty])}{\textup{Var}(q[\tilde{f}^\infty])}
 $$
and it holds that
	$$
\lim_{t\to \infty}\lim_{\varepsilon\to 0} \lambda^*(t,v)=1, \qquad  \lim_{t \to \infty}\lim_{\varepsilon\to 0} \textup{Var}(q^{\lambda^*}[f])=0.
$$

In practice, the optimal value $\lambda^*$ can only be determined numerically. Moreover, accurately computing $\mathbb{E}[q[\tilde{f}]]$ or $\mathbb{E}[q[\tilde{f}^\infty]]$, either exactly or with very small error, is crucial to fully exploit the benefits of the control variate approach.

In order to compute $\textup{Var}(q[\tilde{f})$ and $\textup{Cov}(q[f_\varepsilon],q[\tilde{f}])$ which appear in \eqref{eq:lambda_star} we consider $M$ realizations of $z$ and we use the following unbiased estimators:
\begin{align}
	E_M(q[f_\varepsilon])&=\frac{1}{M}\sum_{k=1}^{M}  q[{f}_\varepsilon^k],\\
	\textup{Var}_M(q[\tilde{f}])&=\frac{1}{M-1}\sum_{k=1}^{M} \left( q[\tilde{f}^k]-\mathbb{E}[\tilde{f}^k]\right)^2, \\ \textup{Cov}_M(q[f_\varepsilon],q[\tilde{f}])&=\frac{1}{M-1}\sum_{k=1}^{M} \left( q[f_\varepsilon^k]-E_M[q[f_\varepsilon]]\right)\left( \tilde{f}^k-\mathbb{E}[\tilde{f}^k]\right).
\end{align}

As far as concern the error introduced by the reconstruction we recall the following Theorem from \cite{pareschi2022mean}.

\begin{theorem}
 The error introduced by the reconstruction function in the control variate approach using DSMC method satisfies
 \begin{align}
 	\lvert\lvert \mathbb{E}[f_\varepsilon](t,\cdot) - E^{\lambda^*}_M [f_{\varepsilon,N,\Delta v}](t,\cdot) \lvert\lvert_{L^p([0,1],L^2(\Omega,L^2([0,1])))}\\
 	\le \frac{\left\| \left(1-\rho^2_{(S,f_\varepsilon),(S,\tilde{f})}\right)^\frac{1}{2}\nu_{(S,f_\varepsilon)}\right\|}{M^\frac{1}{2}}+\frac{\|\sigma_{S,M}\|_{L^p([0,1])}}{N^\frac{1}{2}}+C_{\mathbb{E}[f]}(\Delta v)^q
 \end{align}
 where $\nu^2_{(S,f_\varepsilon)}=\textup{Var}[(S,f)]$ and $\sigma$ is defined in Theorem \ref{thm:rec}.
\end{theorem}

As a result, when the solution of the full model closely approximates that of the control variate, the statistical error caused by uncertainty effectively disappears. This motivates the use of a larger sample size in the state space, aligned with the reconstruction method, to balance the last two error terms in the estimation.

%Moreover, this approach can be extended to a time-dependent control variate method to achieve benefits not only for larger times but also at smaller time scales.
Moreover, instead of using the steady state to improve the Monte Carlo estimate, we employ a time-dependent approximation of the solution $f$ that preserves the moments. 

%The formulation presented above remains valid, however, $f^\infty$ should be replaced with a suitable time-dependent approximation of  $f$ , denoted as $\tilde{f}(t,v,z)$.
\noindent
One possible example of this approach is to consider a BGK-type approximation of the original kinetic equation \cite{HPM2006, herty2020bgk}
$$
\frac{\partial \tilde{f}}{\partial t}=\nu(\tilde{f}^\infty-\tilde{f})
$$
which leads to the following expression for  $\tilde{f}(t,v,z)$:
\begin{equation}\label{eq:lf_BGK}
	\tilde{f}(t,v,z)=e^{-\nu t} f_0(v,z)+(1-e^{-\nu t})\tilde f^\infty(v,z).
\end{equation}

Note that $\tilde f^\infty$, or $\tilde{f}$, described above represents the low fidelity solution whereas the full kinetic model provides the high fidelity solution.

\begin{remark}
This procedure can be extended to multiple multi-scale control variates as shown in \cite{dimarco2020multiscale}. Here we recall briefly the main changes.

We consider $f_1(t,v,z),\dots,f_L(t,v,z)$ approximations of $f(t,v,z)$. We can define the random variable 
\begin{equation}\label{eq:multi}
	f^{\lambda_1,\dots,\lambda_L}(t,v,z)= f(t,v,z) - \sum_{h=1}^L \lambda_h (f_h(t,v,z) - \mathbb{E}[f_h](t,v))
\end{equation}
with $\mathbb{E}[f^{\lambda_1,\dots,\lambda_L}]=\mathbb{E}[f]$ and variance  
\begin{equation}
	\textup{Var}(f^{\lambda_1,\dots,\lambda_L}) = \textup{Var} (f) + \sum_{h=1}^{L} \lambda_h^2 \textup{Var}(f_h) + 2 \sum_{h=1}^{L} \lambda_h \left( \sum_{k=1,k\neq h }^{L}  \lambda_k \textup{Cov}(f_h,f_k) - \textup{Cov}(f,f_h) \right).
\end{equation}

The control variate estimator based on \eqref{eq:multi} reads as:
\begin{equation}
	E_M^\Lambda [f](t,v) = E_M [f](t,v) - \sum_{h=1}^{L} \lambda_h (E_M [f](t,v) - \boldsymbol{f}_h (t,v) )
\end{equation}
where $\Lambda=(\lambda_1,\dots,\lambda_L)^T$, $\boldsymbol{f}_h (t,v)$ is an approximation of $\mathbb{E}[f_h](t,v)$.

Note that to have a more stable estimator, the surrogates models can be chosen with an increasing level of fidelity. Under this assumption the control variate $f_1$ represents the less accurate model while $f_L$ is the closer to the full model $f$.
In particular, we estimate $\mathbb{E}[f]$ with $M_L$ samples using $f_L$ as control variate:
\begin{equation}
	\mathbb{E}[f] \approx E_{M_L} [f] - \hat{\lambda}_L (E_{M_L}[f_L] - \mathbb{E}[f_L]).
\end{equation}
Then, we use $M_{L-1}\gg M_L$ samples to estimate $\mathbb{E}[f_L]$ and we consider $f_{L_1}$ as control variate:
\begin{equation}
\mathbb{E}[f_L] \approx E_{M_{L-1}} [f_L] - \hat{\lambda}_{L-1} (E_{M_{L-1}}[f_{L-1}] - \mathbb{E}[f_{L-1}]).	
\end{equation}

In the same way, we can recursively write the remaining expectations of the control variates using respectively $M_{L-3} \ll M_{L-4} \ll \dots \ll M_1$ until the final estimate: $\mathbb{E}[f_1]\approx E_{M_0}[f_1]$ with $M_0 \gg M_1$.

Summarizing the estimate:
\begin{align}\label{eq:hier}
E_L^{\hat{\Lambda}} [f] = & E_{M_L} [f_{L+1}] - \sum_{h=1}^{L} \lambda_h 
(E_{M_h} [f_h] - E_{M_{h-1}}[f_h])\\
= &  \lambda_1 E_{M_0} [f_{1}] + \sum_{h=1}^{L} (\lambda_{h+1} 
E_{M_h} [f_{h+1}] - \lambda_h E_{M_{h}}[f_h])
\end{align}
where we consider $\lambda_{h}=\prod_{j=h}^{L} \hat{\lambda}_j$, for $h=1,\dots,L$ and $\lambda_{L+1}=1.$

Further details can be found in \cite{dimarco2020multiscale}. 
\end{remark}

\subsection{Multi-level approach}\label{sec:multilevel}

In order to reduce the variance related to the approximation, another technique is the Multilevel Monte Carlo (MLMC) \cite{giles2006improved,giles2015multilevel,heinrich2001multilevel} technique which relies on approximating the high fidelity simulations with a hierarchy of models at varying levels of resolution. In particular, it evaluates most realizations on a coarse level of accuracy, with only a few realizations on a high level of solution refinement.

The general form of the MLMC approach is:
\begin{equation}\label{eq:MLMC}
	\mathbb{E}[f]= \mathbb{E}[f_0] + \sum_{h=1}^{L} \mathbb{E}[f_h-f_{h-1}]
\end{equation}

where $h$ denotes the refinement level, $f_h$ indicates the reconstruction of $f$ made at level $h$. In our case, $h$ stands for the  samples used at each level in the stochastic space.

The MLMC method works if the variances $\textup{Var}(f_h-f_{h-1}) \to 0$ as $h\to\infty$, which occurs when $f_h$ and $f_{h-1}$ approximate the same random variable $f$.

Here we want to highlight the analogies between the multi-fidelity approach and the MLMC by employing a hierarchy of discretizations of the kinetic equation.

For example, on a cartesian grid this approach aims to construct a sequence of velocity discretizations with corresponding mesh width $\Delta v_h$ such that $\Delta v_h = 2^{1-h}(\Delta v_1)$ for $h=1,\dots,L$, where $\Delta v_1$ is the mesh width at the coarsest resolution.
Let us denote by $f_h(t,v,z)$ the continuous representation of the corresponding numerical solution at time $t$ obtained with the deterministic method using mesh width $\Delta v_h$. 
Under these assumptions, we recover the classical MLMC estimator \eqref{eq:MLMC} by setting $\lambda_h=1$:
\begin{equation}\label{eq:MLest}
	E^\mathbf{1}_L[f](t,v) = E_{M_0}[f_1] + \sum_{h=1}^{L} (E_{M_h}[f_{h+1} - f_h] )
\end{equation}
where $\mathbf{1}=(1,\dots,1)^T$.
 
Note that, exploiting the results from the multi-fidelity approach, we can construct a optimal (quasi-optimal) version of the MLMC just using the optimal (quasi-optimal) values for $\lambda_h$, $h=1,\dots,L$.

Moreover, as in the MLMC technique, the multi-fidelity estimator  where the surrogates models are chosen hierarchically \eqref{eq:hier}, requires the  largest number of samples, $M_0$, on the less accurate model $f_1$, where samples are cheaper, while only a small number $M_L$ of samples is needed on the full model. 

The two main differences between recursive multi-fidelity estimator \eqref{eq:hier} and the MLMC approach \eqref{eq:MLest} are the use  of low-fidelity models as control variates instead of a hierarchy of discretizations, and the use of quasi-optimal values intead of setting $\lambda_h =1$ for $h=1,\dots,L$.

\section{Numerical results}\label{sec:appl}

In this section, we aim to combine the stochastic kinetic traffic flow model with the control variate approaches, presenting and discussing the numerical results. %Different low fidelity models will be discussed and compared. 
We denote by $N$ the number of samples in the physical space $v\in [0,1]$ and by $M$ the number of samples in the stochastic space, $z_1,\dots,z_M$. The discretization points for $v$ in the histogram reconstruction are set to $N_v=10^2$. Additionally, we assume the uncertainty, which lies in the interaction kernel, to be uniformly distributed, with $z \sim \mathcal{U}(1,3)$ as in \cite{tosin2019kinetic}.

The high-fidelity model is numerically solved using the Direct Simulation Monte Carlo method, considering the interaction rules described in \eqref{eq:update1}-\eqref{eq:update2}. 
The value of the density in \eqref{eq:interaction} is set to $\rho=0.4$. This choice is motivated by the fact that this value typically represents the switching point between the free flow regime and congested traffic, also known as the critical density.

In the following, the term \textit{reference solution} refers to the accurate solution computed using the high-fidelity model using a large number of samples $N_r=10^6$ in the physical space. The reference solution is computed offline, once and for all. In particular, we used a stochastic collocation approach with 10 nodes to compute the expectation of the reference solution in time, employing the Gauss-Legendre quadrature rule. We denote by $E_{GL}[f_{ref}]$ the solution computed in this way.  

To compare the different approaches, we will focus on the relative error, which is computed as 
\begin{equation}\label{eq:err}
\frac{||E_{GL}[f_{ref}]-E_M[f_N]||_2}{||E_{GL}[f_{ref}]||_2},
\end{equation}
 where $f_N$ is the histogram reconstruction of the solution of the kinetic equation, $f_{ref}$ is the histogram reconstruction of the reference solution and $E_M$ are the different Monte Carlo estimators used to evaluate the expectation with respect to the uncertainty. 
Moreover, we consider two different regimes, setting $\varepsilon=1$ and $\varepsilon = 0.003$. In the second case, the surrogate model is expected to perform better, as its steady state aligns more closely with the steady state to which the full kinetic system converges.

\subsection{Two and three levels Monte Carlo}

Numerically, we consider two and three levels Monte Carlo, i.e. with $L=2$ and $L=3$ in \eqref{eq:MLMC} respectively. We set $N=10^4$ particles in the physical space and $M=30$ in the stochastic space. In particular, for each level we refine the stochastic space by doubling the sampling while taking half of the samples in the velocity space. Thus, for the two-level Monte Carlo we consider 
\begin{equation}
	E_M[f_N]= E_M[f_N] - \left[ E_{M}[f_{\nicefrac{N}{2}}] - E_{2M}[f_{\nicefrac{N}{2}}] \right].
\end{equation} 
In the three levels case, we also consider $E_{4M}[f_{\nicefrac{N}{4}}]$. In particular we replace $E_{2M}[f_{\nicefrac{N}{2}}]$ with
\begin{equation}
	E_{2M}[f_{\nicefrac{N}{2}}] = E_{2M}[f_{\nicefrac{N}{2}}] - \left[ E_{2M}[f_{\nicefrac{N}{4}}] - E_{4M}[f_{\nicefrac{N}{4}}] \right].
\end{equation}

\begin{figure}[h!]
	\centering
	%\animategraphics[width=0.48\textwidth]{5}{fig/Frames_video/FIG}{1}{40}\quad
	%\includegraphics[width=0.48\textwidth]{fig/T0_MLMC.pdf}
	\includegraphics[width=0.48\textwidth]{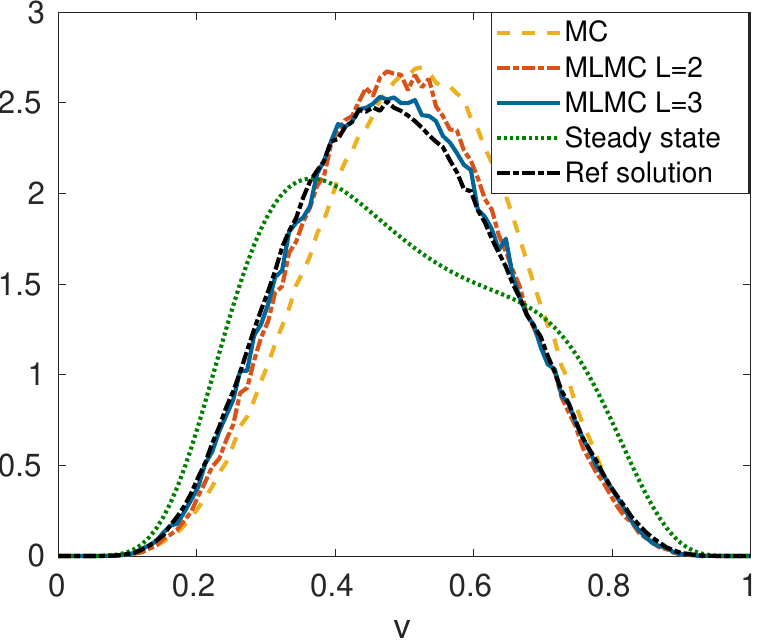}
	\includegraphics[width=0.48\textwidth]{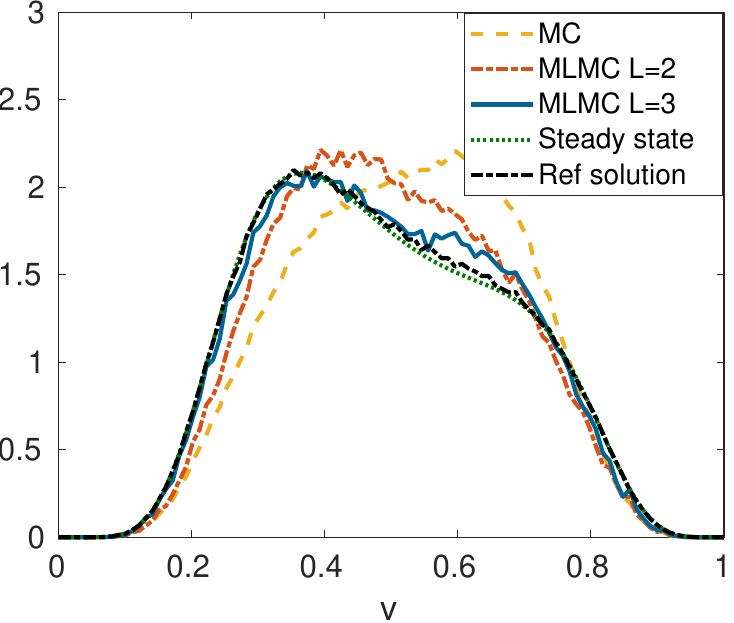}
	\caption{{\bf MLMC approach}. Expected value of the kinetic density $\mathbb{E}[f]$ for $\varepsilon=0.003$ at different times. Solution computed by Monte Carlo (yellow dashed line), MLMC with $L=2$ (red dashed dotted line), MLMC with $L=3$ (blue line), steady state (green dotted line) and reference solution (black dashed dotted line), at times $T=\frac{1}{4}T_f$ (left) and $T=T_f$ (right), for $z \sim \mathcal{U}([1,3])$. }
	\label{fig:sol_T_MLMC}
\end{figure}

\begin{figure}[h!]
	\centering
	%\animategraphics[width=0.48\textwidth]{5}{fig/Frames_video/FIG}{1}{40}\quad
		\includegraphics[width=0.505\textwidth]{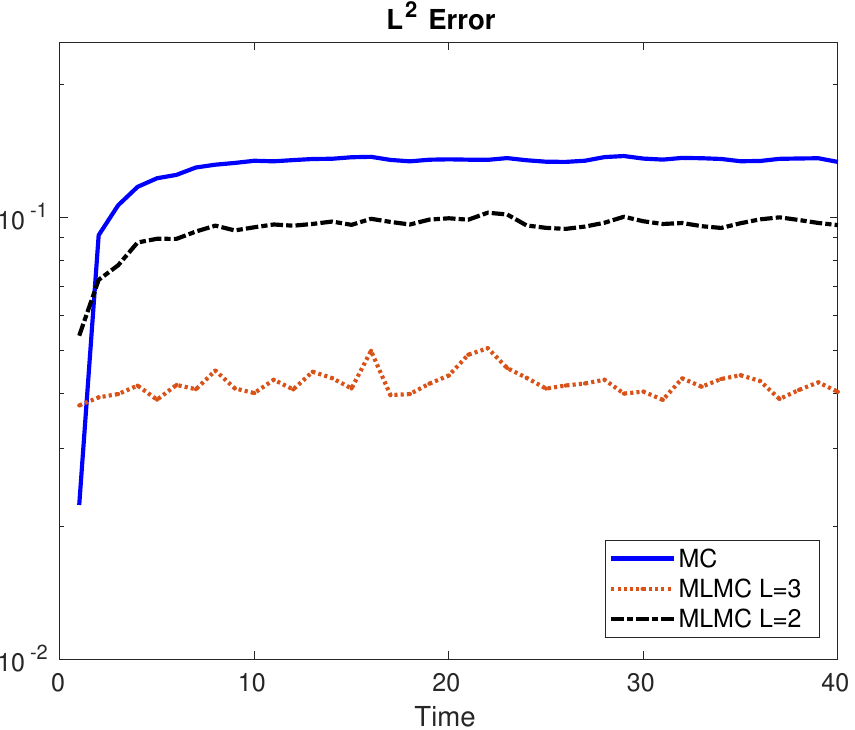}
		\includegraphics[width=0.473\textwidth]{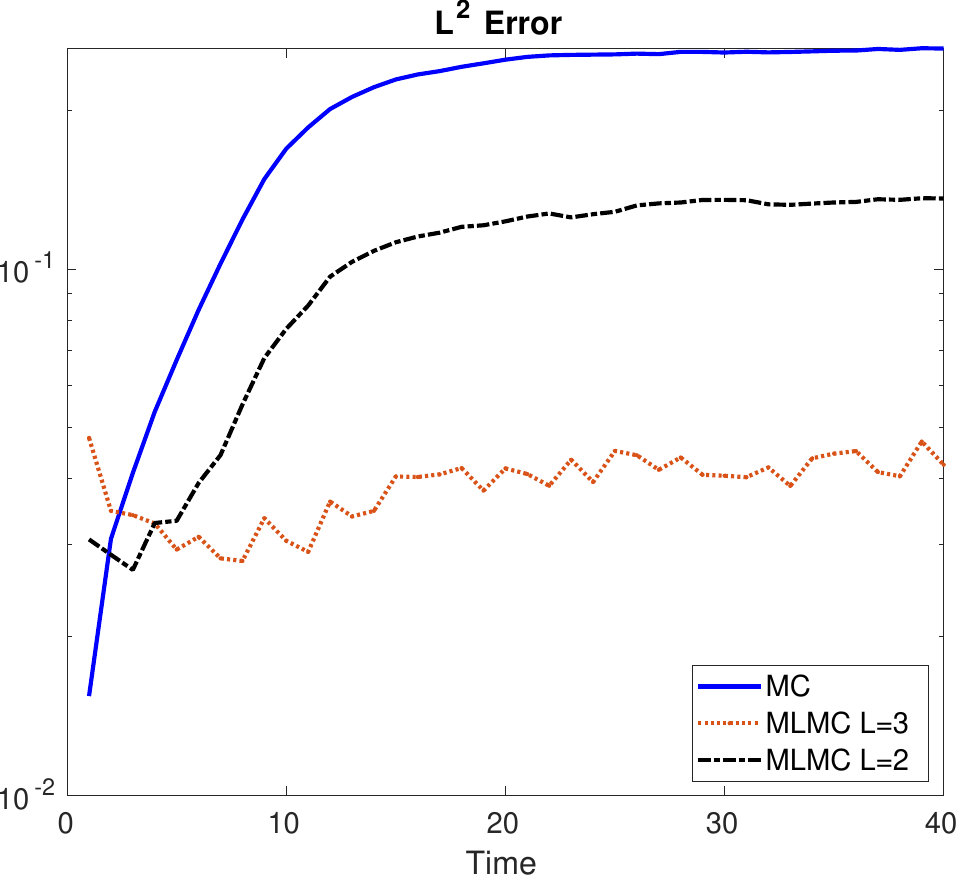}
	\caption{{\bf MLMC approach}. Comparison between the relative errors, computed following \eqref{eq:err}, of the Monte Carlo method (blue line), MLMC with $L=2$ (black dashed line) and MLMC with $L=3$ (red dotted line) for $\varepsilon=1$ (left) and $\varepsilon=0.003$ (right).}
	\label{fig:comp_MLMC}
\end{figure}

In Fig. \ref{fig:sol_T_MLMC} we compare the expectation of the kinetic density $\mathbb{E}[f]$ for $\varepsilon=0.003$ obtained using different Monte Carlo approaches: the standard MC, the MLMC with $L=2$ and $L=3$. Starting from a uniform distribution at time $T=0$ we observe at a final time $T_f=40$ (right) that the MLMC with $L=3$ solution shows significant improvement, closely approximating the reference solution compared to the two-level and classic MC solutions.
This behaviour is even more evident when we focus on the comparison between the relative errors, as illustrated in Fig. \ref{fig:comp_MLMC}.  We compute the $L^2$-norm of the difference between the means of the quantity of interest and the reference solution, normalized by the norm of the mean of the reference solution, as stated in eq. \eqref{eq:err}.
In particular, we show the relative errors of the solutions computed by the standard Monte Carlo method, the MLMC for $L=2$ and $L=3$.  Results are presented for two scenarios: $\varepsilon=1$ (left), corresponding to the kinetic regime and $\varepsilon=0.003$ (right), representing a regime near the mean-field limit. As expected, the MLMC with three levels performs better than the two-level one. Moreover, the order of the relative error is qualitatively the same for both choices of $\varepsilon$.

\subsection{Bi-fidelity method with low fidelity approximated steady state}

Let us consider the steady state \eqref{eq:steady} as low fidelity model in \eqref{eq:cv_steady}. This implies that the cost of performing the low-fidelity model over the whole time range is reduced to a single function evaluation.
To assess the performance of the different control variate approaches, we analyze the evolution of the expected velocity distribution over time, as given by \eqref{eq:cv_steady} and the relative error over time. We use $M=20$ and perform $10^4$ evaluations of the low fidelity to compute $\boldsymbol{f}^\infty = \mathbb{E}[{f}^\infty]$. 

\begin{figure}[h!]
	\centering
	%\animategraphics[width=0.48\textwidth]{5}{fig/Frames_video/FIG}{1}{40}\quad
	%\includegraphics[width=0.48\textwidth]{fig/T0_SS.pdf}
	\includegraphics[width=0.48\textwidth]{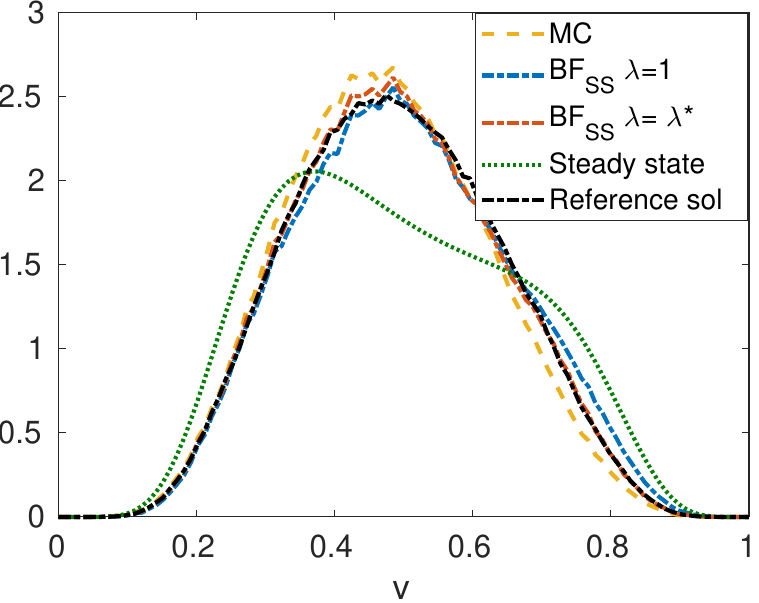}
	\includegraphics[width=0.48\textwidth]{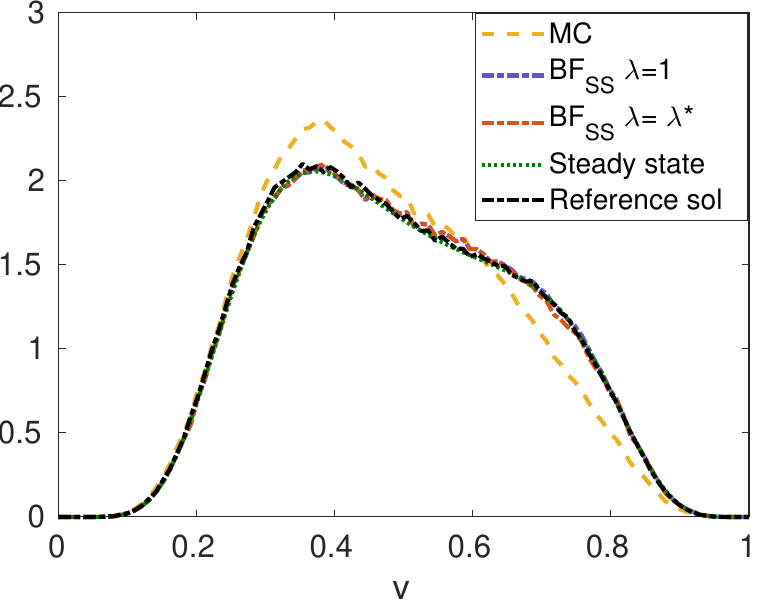}
	\caption{{\bf Multifidelity: steady state}.
	Expected value of the kinetic density $\mathbb{E}[f]$ for $\varepsilon=0.003$ at different times. Solution computed by Monte Carlo (yellow dashed line), Bi Fidelity (BF) approach with low fidelity the steady state and $\lambda=1$ (blue dashed dotted line) $\lambda=\lambda^*$ (red dashed dotted line), steady state (green dotted line) and reference solution (black dashed dotted line), at times $T=\frac{1}{4}T_f$ (left) and $T=T_f$ (right), for $z \sim \mathcal{U}([1,3])$.}
	\label{fig:sol_T}
\end{figure}

In Fig. \ref{fig:sol_T}, we compare the solutions obtained from the bi-fidelity approaches with $\lambda=1$ and $\lambda=\lambda^*$ the estimated optimal value. It is evident that over time, the bi-fidelity approach consistently remains closer to the reference solution compared to the Monte Carlo approach, and convergence to the steady state is achieved in the bottom-right plot, while the Monte Carlo method lags behind. A slightly better performance of the method with the optimal value $\lambda^*$ is observed for short times at $T=\frac{1}{4}T_f$.

This improvement in performance is also reflected in the evolution of the relative error over time. We compute the $L^2$-norm of the difference between the means of the quantity of interest and the reference solution, normalized by the norm of the mean of the reference solution, as stated in eq. \eqref{eq:err}.
In Fig. \ref{fig:comparison} we can see the comparison between two choices of the parameter $\lambda$: $\lambda=1$ and the optimal value $\lambda=\lambda^*$. 
The solution computed by the bi-fidelity approach with low fidelity the steady state with $\lambda^*$ performs better with respect to the solution obtained with $\lambda=1$, which gives good results only at later times. This is expected, as theory
suggests that $\lambda^*$ approaches $1$ as $t$ increases. However, in a regime far from the Fokker-Planck limit, such as $\varepsilon = 1$, this approach performs significantly worse, as shown in Fig. \ref{fig:comparison} (left).

 %It is interesting to study how this limiting process works. 
%Indeed, in Fig.\ref{fig:rate_T} (right) we plot the evolution of  the correlation coefficient $\rho$ over time with respect to $v$.  The color map, from blue to yellow, indicates the value of $\rho(t,v)$.  We observe a steep increase from zero to values approaching 1 within a short time interval, after which the value adjusts and stabilizes near 1. This implies the convergence of the $\textup{Var}(f^{\lambda^*})$ to $0$ and of $\lambda^*$ to $1$, as stated in the Theorem \ref{thm:multi}.

%\begin{figure}[ht!]
%	\centering
%	%\animategraphics[width=0.48\textwidth]{5}{fig/Frames_video/FIG}{1}{40}\quad
%	\includegraphics[width=0.46\textwidth]{fig/error_SS.pdf}
%%	\includegraphics[width=0.48\textwidth]{fig/Lambda_star_Tempi.pdf}
%%	\includegraphics[width=0.515\textwidth]{fig/rho.pdf}
%	\caption{Comparison between the relative errors, computed following \eqref{eq:err}, of the Monte Carlo method (blue line) and the Bi Fidelity approach with low fidelity the steady state and $\lambda=1$ (green dashed line) and $\lambda=\lambda^*$ (red dotted line).  }
%	\label{fig:rate_T}
%\end{figure}

\begin{figure}[ht!]
\centering
%\animategraphics[width=0.48\textwidth]{5}{fig/Frames_video/FIG}{1}{40}\quad
\includegraphics[width=0.49\textwidth]{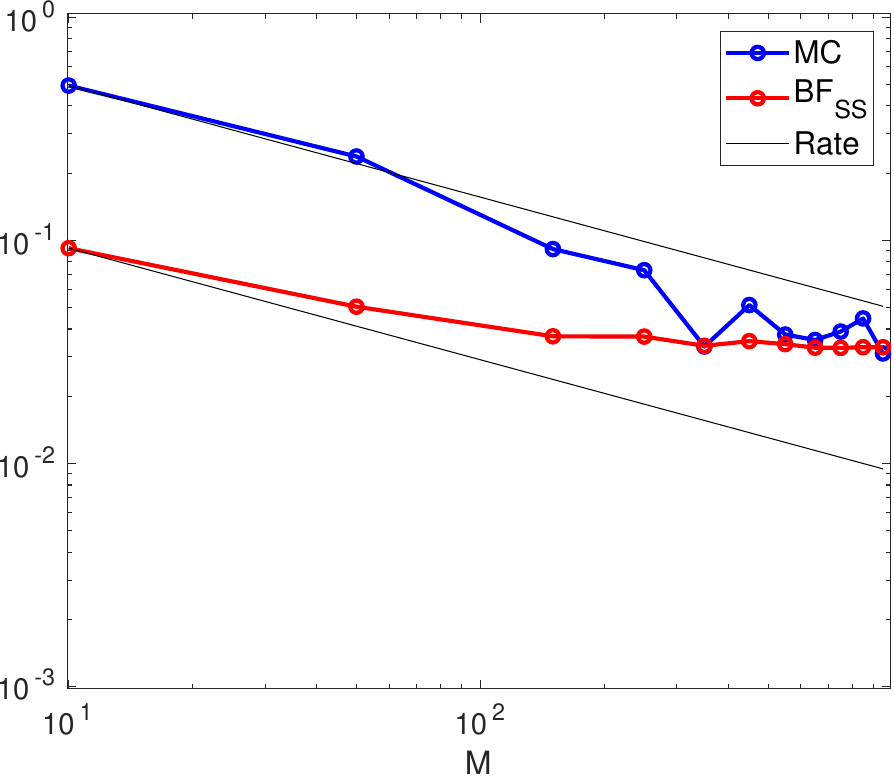}
\includegraphics[width=0.48\textwidth]{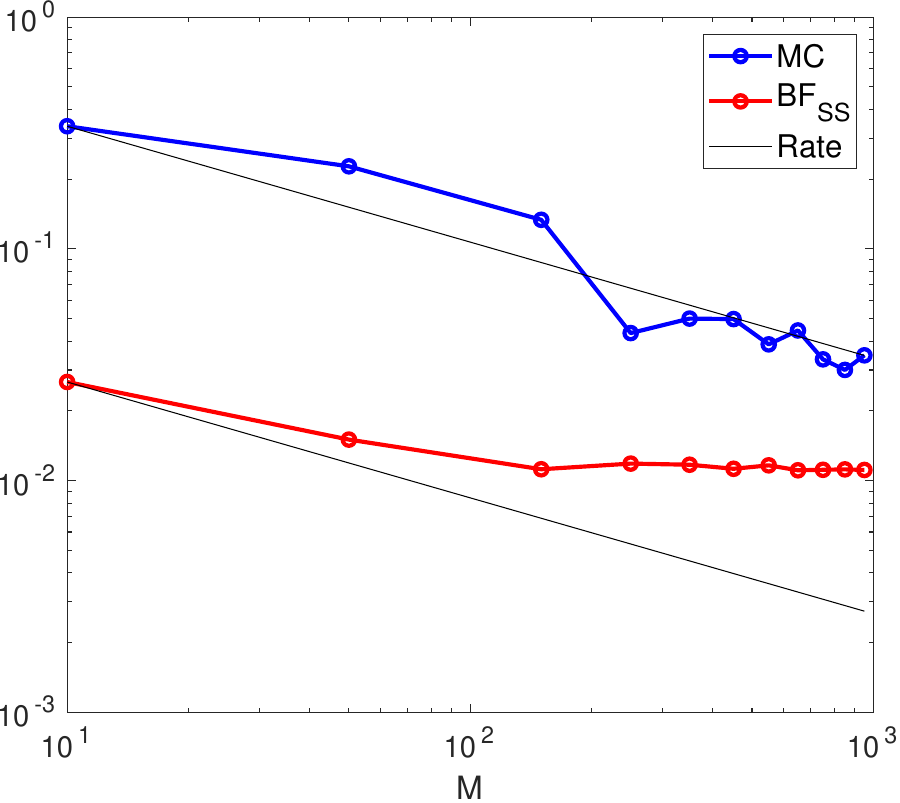}
\caption{{\bf Multifidelity: steady state}. Convergence rate for $\varepsilon=0.003$ of the Monte Carlo method (blue line) and the Bi Fidelity approach (red line) with respect to the stochastic number of samples $M$ for a fixed number of samples in the physical space $N=10^4$ (left), $N=10^5$ (right). }
\label{fig:rate_M}
\end{figure}

It is also interesting to investigate the behavior of the error concerning the number of samples in the stochastic space. We fix the time $t=T_f$ and compute the error for various value of $M$. In Fig. \ref{fig:rate_M}, we set $N=10^4$ (left) and $N=10^5$ (right). One can observe that while the error of the  Monte Carlo method decreases as $M^{\frac{1}{2}}$, in the bi-fidelity approach the error reaches a plateau, which is caused by the interaction between the nested Monte Carlo methods, one running on the physical space and the one on the stochastic space. 
This can be observe since for different values of $N$ the qualitative plot remains the same, although the level of the plateau changes. 
Analyzing the effects of nested Monte Carlo methods will be the focus of future work.

\subsection{Bi fidelity method with low-fidelity BGK approximation}

\begin{figure}[ht!]
	\centering
	%\animategraphics[width=0.48\textwidth]{5}{fig/Frames_video/FIG}{1}{40}\quad
	\includegraphics[width=0.48\textwidth]{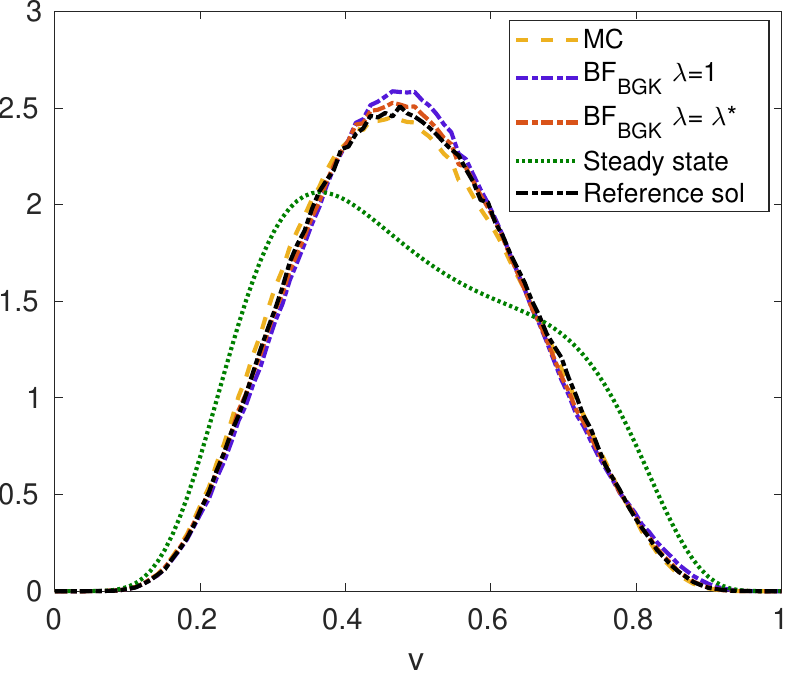}
	\includegraphics[width=0.49\textwidth]{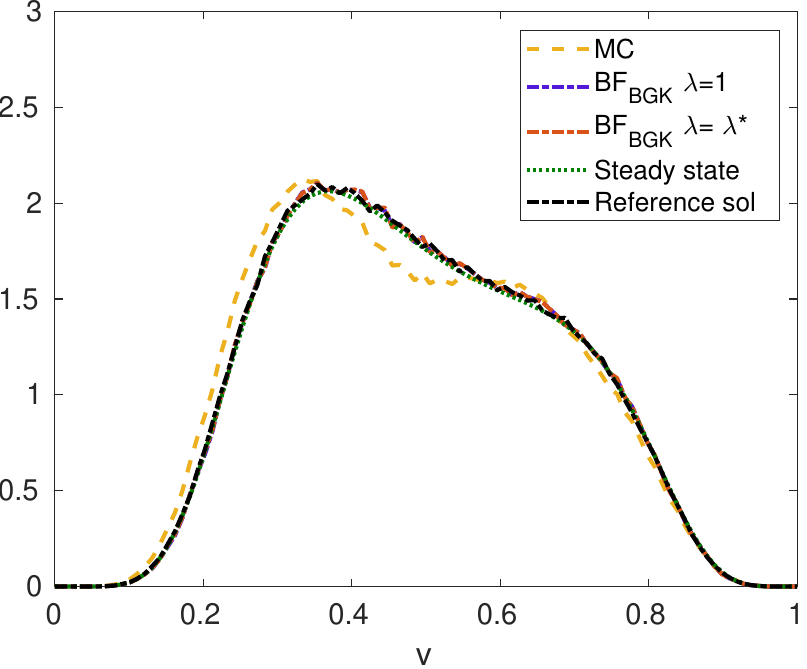}
	\caption{{\bf Multifidelity: BGK model}. Expected value of the kinetic density $\mathbb{E}[f]$ for $\varepsilon=0.003$ at different times. Solution computed by Monte Carlo (yellow dashed line), Bi Fidelity (BF) approach with low fidelity BGL approximation and $\lambda=1$ (blue dashed dotted line) $\lambda=\lambda^*$ (red dashed dotted line), steady state (green dotted line) and reference solution (black dashed dotted line), at times $T=\frac{1}{4}T_f$ (left) and $T=T_f$ (right), for $z \sim \mathcal{U}([1,3])$. }
	\label{fig:sol_T_BGK}
\end{figure}

\begin{figure}[ht!]
	\centering
	%\animategraphics[width=0.48\textwidth]{5}{fig/Frames_video/FIG}{1}{40}\quad
	\includegraphics[width=0.485\textwidth]{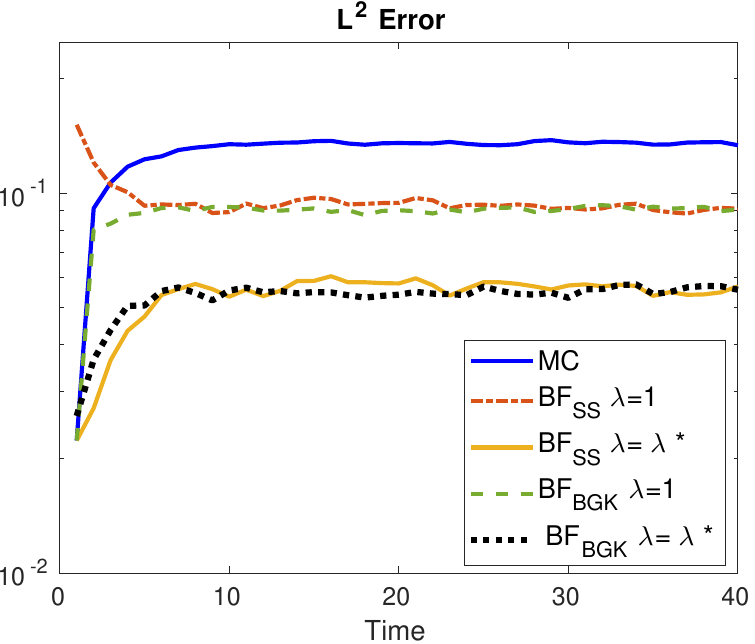}
	\includegraphics[width=0.5\textwidth]{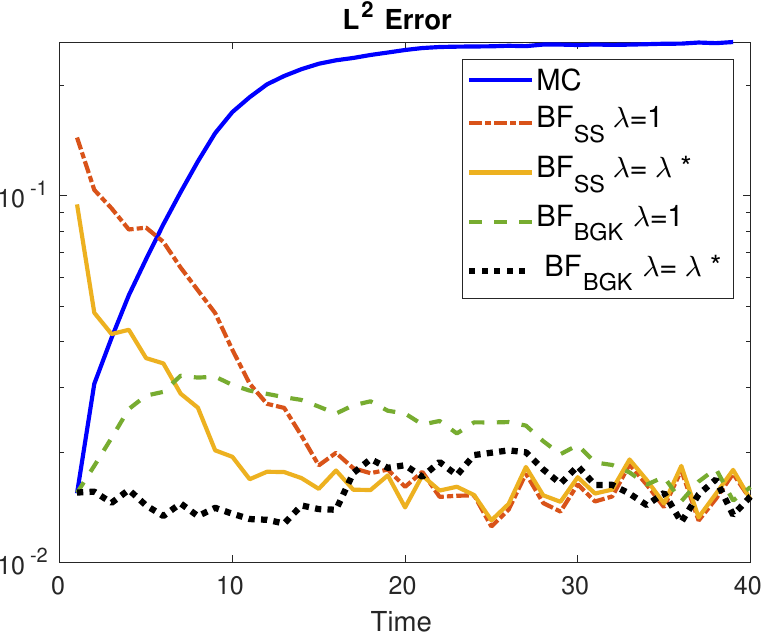}%comp_BF_eps_piccolo
	\caption{{\bf Multifildelity approaches.} Comparison between the relative errors, computed following \eqref{eq:err}, of the Monte Carlo method (blue line) and the Bi-Fidelity approach with steady state low fidelity, with both $\lambda=1$ (red dashed-dotted line) and $\lambda=\lambda^*$ (yellow line), and BGK type low fidelity with both $\lambda=1$ (green dashed line) and $\lambda=\lambda^*$ (black dotted line) for $\varepsilon=1$ (left) and $\varepsilon=0.003$ (right). }
	\label{fig:comparison}
\end{figure}

As low fidelity model, we now consider \eqref{eq:lf_BGK}, which incorporates information from both the initial data and the steady state. In this case the choice of the relaxation parameter $\nu$ plays a major rule and in principle can be optimized to maximize correlations. Here we simply assume $\nu=\frac{1}{2}$. 
In Fig.\ref{fig:sol_T_BGK}, we compare the solutions obtained with the control variate approach using different values of $\lambda$, namely $\lambda=1$ and the estimated optimal value $\lambda=\lambda^*$. The number of samples is $M=30$.
Similar to Fig. \ref{fig:sol_T}, the MC solution remains significantly distant from the reference solution, even when the steady state is achieved, due to the few samples used. We note that for $T=\frac{1}{4}T_f$ the solution computed with $\lambda^*$ is more accurate than the one computed with $\lambda=1$ although this difference diminishes over time. 
This is even more evident in Fig. \ref{fig:comparison} (left), where the difference between the error of the solution computed with $\lambda=1$ is very large at the beginning, since the low fidelity heavily weights the initial data. However, this difference vanishes over time due to the convergence of $\lambda$ and the decreasing influence of the initial data on the low-fidelity model. 
Note that the uncertainty lies in the interaction kernel, which allows the BGK to benefit from knowledge of the initial data at smaller time steps, but this advantage diminishes over longer times.

%\begin{figure}[ht!]
%	\centering
%	%\animategraphics[width=0.48\textwidth]{5}{fig/Frames_video/FIG}{1}{40}\quad
%	\includegraphics[width=0.46\textwidth]{fig/errore_tempi_bgk.pdf}
%%	\includegraphics[scale=0.47]{fig/rho_BGK.pdf}
%	\caption{Comparison between the relative errors, computed following \eqref{eq:err}, of the Monte Carlo method (blue line) and the Bi-Fidelity approach with BGK type low fidelity and $\lambda=1$ (green line) and $\lambda=\lambda^*$ (red line). %The plot on the right shows the evolution of the correlation coefficient, $\rho(t,v)$. 
%	 }
%	\label{fig:rate_BGK}
%\end{figure}

%The evolution of $\rho^*(v,t)$ is quite similar to the steady-state case, see Fig. \ref{fig:rate_T}(right), though it tends to reach 1 during the very early stages, as one can see in Fig. \ref{fig:rate_BGK} (right).

\begin{figure}[ht!]
	\centering
	%\animategraphics[width=0.48\textwidth]{5}{fig/Frames_video/FIG}{1}{40}\quad
	\includegraphics[width=0.495\textwidth]{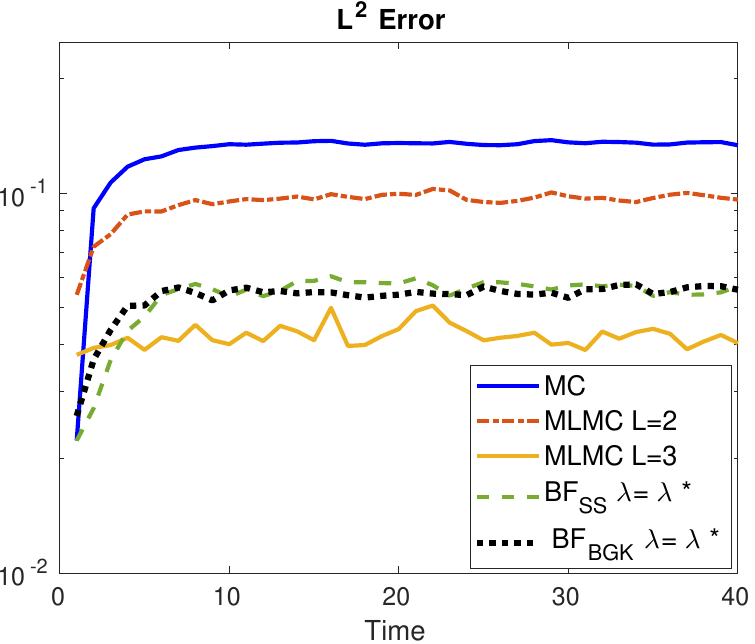}
	\includegraphics[width=0.495\textwidth]{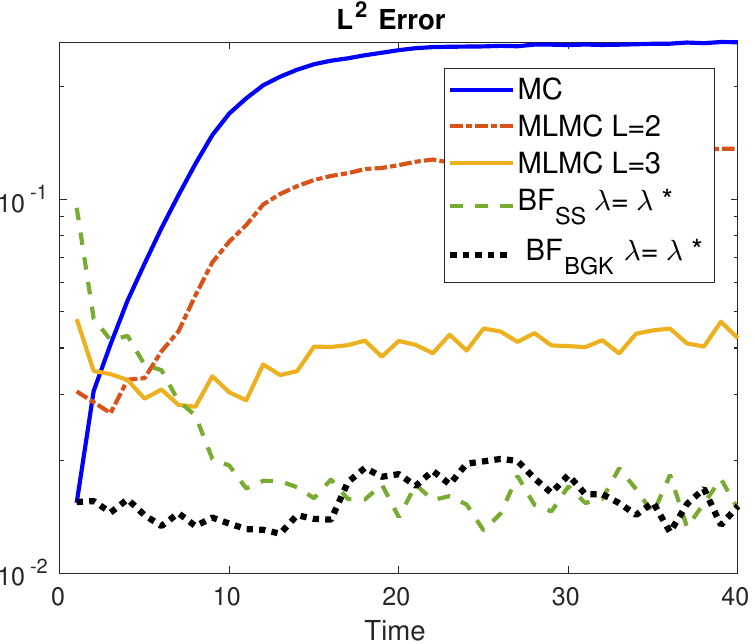}%comp_eps_picc
	\caption{{\bf Multifildelity vs MLMC.} Left: comparison between the relative errors, computed following \eqref{eq:err}, of the Monte Carlo method (blue line) and the Multi-level Monte Carlo approach with $L=2$ (red dashed-dotted line) and $L=3$ (yellow line) and Bi-Fidelity approach with steady state low fidelity (green dashed line) and BGK type low fidelity (black dotted line) both with $\lambda=\lambda^*$, for $\varepsilon=1$ (left) and $\varepsilon=0.003$ (right).}
	\label{fig:comparison_fin}
\end{figure}

Moreover in Fig. \ref{fig:comparison} we can see the comparison  between the solution obtained with the control variate approach using the  steady state low fidelity and BGK type low fidelity, both with $\lambda=1$ and $\lambda=\lambda^*$. 
We note that for longer times the bi-fidelity with the BGK low fidelity performs better than considering as surrogate model only the steady state. 
As expected, the solutions obtained by considering as low fidelity model the steady state and the one obtained by the BGK low fidelity converge to the same result, since the model itself reaches the steady state and the optimal $\lambda$ is 1. 
In Fig. \ref{fig:comparison} (right), we compare the relative errors of the Monte Carlo method with those of the bi-fidelity approach across different values of $\varepsilon$. Specifically, we examine the cases $\varepsilon = 1$, corresponding to a regime closer to the kinetic scale and $\varepsilon = 0.003$, representing a regime near the Fokker–Planck scale. 
%The analysis focuses on two scenarios: $\lambda = 1$ for $\varepsilon=0.003$ and green dashed line for $\varepsilon=1$) and $\lambda = \lambda^*$ (yellow solid line for $\varepsilon=0.003$ and black dotted line for $\varepsilon=1$). \
As expected, for longer times, the approximation using the steady-state model performs better for smaller values of $\varepsilon$, as the system converges toward the steady state employed as the surrogate model. In contrast, for $\varepsilon = 1$, there is less correlation between the steady state and the kinetic regime, leading to reduced accuracy. Finally, in Fig. \ref{fig:comparison_fin} we summarize the results of the MLMC and Bi-fidelity approaches. Overall, for $\varepsilon=0.003$, there is an improvement in accuracy of Bi-fidelity over over MLMC. In particular, when $\lambda$ is optimal and  the BGK model serves as the low-fidelity surrogate, the results are further enhanced (see Fig. \ref{fig:comparison_fin}, right). On the other hand, for $\varepsilon=1$ (Fig. \ref{fig:comparison_fin}, left), the MLMC with $L=3$ (yellow line) shows better performance. This highlights the importance of selecting an appropriate surrogate model to fully exploits the benefits of the multi-fidelity approach.

\section{Conclusions}
This work presents a comprehensive framework for employing control variate strategies in traffic flow models influenced by uncertain factors. Our focus was on designing effective approaches to deal with uncertainties in vehicles interactions, therefore we focused on a homogeneous kinetic model of traffic flow, through which we analyzed and compared the performance of multi-fidelity and multi-level Monte Carlo methods. 
By examining their structural similarities and operational distinctions, we introduced both approaches in the case when a Monte Carlo method 
is adopted in the combined physical-random space. Numerical results indicate that, in cases where low-fidelity modeling surrogates are available, like for example approximate steady state solutions or simplified BGK interaction operators, the multi-fidelity approach may provide a more favorable error relative to the multilevel Monte Carlo method. This advantage in error reduction demonstrates the potential of multi-fidelity methods to improve the efficiency and accuracy of uncertainty quantification in traffic flow simulations. As expected, in regimes where the low fidelity model loses correlations with the high fidelity one, a MLMC approach may yield better results.

Upcoming research will expand this unified framework to tackle non - homogeneous models, which introduce additional spatial and temporal complexities. Addressing these challenges will require enhanced high-fidelity and low-fidelity modeling techniques to accurately capture the variation in flow properties across different traffic conditions. The design of hybrid methods which combine both techniques will be also the subject of further studies. Another key direction will involve a more in-depth analysis of the synergy between the two Monte Carlo methods — one dedicated to the kinetic equation itself and the other focused on managing uncertainties. Understanding this interplay will be crucial for refining the efficiency of these approaches, potentially leading to even more precise and computationally efficient solutions for uncertainty quantification in complex traffic models.

%------
% Insert acknowledgments and information
% regarding funding at the end of the last
% section, i.e., right before the bibliography.
%------

%\begin{ack}
%We thank Federica Ferrarese for the invaluable support.
%\end{ack}

\vspace{2cm}
\textbf{Fundings}:
This work has been written within the activities of GNCS group of INdAM (Italian National Institute of High Mathematics) with Project \\ CUP\_E53C23001670001. The research of LP has been supported by the Royal Society under the Wolfson Fellowship ``Uncertainty quantification, data-driven simulations and learning of multiscale complex systems governed by PDEs". The partial support by ICSC -- Centro Nazionale di Ricerca in High Performance Computing, Big Data and Quantum Computing, funded by European Union -- NextGenerationEU and by MIUR-PRIN Project 2022, No. 2022KKJP4X, ``Advanced numerical methods for time dependent parametric partial differential equations with applications" is also acknowledged.

%------
% Insert the bibliography.
%------

\bibliographystyle{siam}
\bibliography{ref}

\end{document}